\documentclass{article}
\usepackage[utf8]{inputenc}

\newtheorem{theorem}{Theorem}[section]

\newtheorem{conjecture}[theorem]{Conjecture}

\newtheorem{definition}[theorem]{Definition}

\newcommand{\cf}{\mathcal{F}}
\newcommand{\F}{\mathcal{F}}

\title{Strongly maximal matchings and strongly minimal covers}
\author{Ron Aharoni }
\date{June 2022}

\begin{document}
\maketitle
\begin{abstract}
    This is a not-to-be-journal-published paper, aimed to serve as reference. It is a summary of the main ideas on the topic appearing in the title, and an opportunity to state  correctly the main conjecture in the field.
\end{abstract}

\section{Motivation}
In finite combinatorics, duality results are formulated as ``$\max (LP)= \min (D)$'', where (LP) and (D) are dual linear programs. In the  infinite case ``maximal size'' and ``minimal size'' (where ``size'' is cardinality) are weak notions, rendering the min-max equality trivial. Erd\H{o}s realized the ``right'' formulation, in his conjectured infinite 
version of Menger's theorem, that has since been proved \cite{ab}. In the language of Linear Programming, it is the existence of a pair (solution of (LP), solution of (D))  satisfying the complementary slackness conditions. 

\begin{theorem}\label{ab}
For any two sets of vertices, $A,B$, in a possibly infinite directed graph,   there exist a family $\F$ of vertex-disjoint $A-B$ paths and an $A-B$-separating set $S$ of vertices, such that 

(a)~ Every path in $\F$ contains precisely one vertex from $S$, and 

(b)~Every vertex in $S$ lies on some path from $\F$.
\end{theorem}

One proof in the finite case 
starts from a  family of  disjoint $A-B$ paths of maximal size, and uses the maximality to generate the separating set $S$. In the infinite case, a stronger version of maximality is needed:

\begin{definition}
Let $H$ be a  hypergraph. A set $K \in H$ is called {\em strongly maximal} (resp. {\em strongly minimal}) if for every $L \in H$ we have $|L \setminus K| \le |K\setminus L|$ (resp.  $|K \setminus L| \le |L\setminus K|$).
\end{definition}

A standard alternating paths argument shows that in Theorem \ref{ab} every strongly maximal set of disjoint $A-B$ paths $\cf$ has an $A-B$ separating set $S$ as in the theorem. On each path $F \in \cf$ you choose the last vertex lying on an alternating path starting at $A \setminus V(\cf)$, and the first vertex of $F$ if there is no such alternating path meeting $F$. The problem is that there is no simple way of constructing such a strongly maximal set. It just emerges from the proof of the theorem. 

Interestingly, it is not possible to start from the other direction, the strongly minimal separating set. 
Take the  bipartite graph with side $A$ being (a copy of) $\omega$
and side $B =\omega \setminus \{0\}$. 
Connect $i$ in $A$ to $i$ in $B$
except for $0$, that is connected to all of $B$.
Let $S =A$. 

To make the conjecture true, we have to  replace the condition that $|T\setminus S| \ge |S\setminus T|$ by the more refined:

``there is no $A-B$-separating set $S$ such that 
$T \setminus S$ is linkable in the graph strictly into $S\setminus T$''.
(``linkability'' is by a set of vertex-disjoint paths).  
\section{ The main conjecture on strong maximality and minimality, and the fish bone conjecture}

A natural place where the notions of strong maximality and strong minimality appear is matchings and covers. A {\em matching} in a hypergraph $H$ is a set of disjoint edges (sets). The notion of  ``cover'' is used in two senses: a {\em vertex cover} is a set of vertices intersecting every edge of the hypergraph, and an {\em edge cover} is a set of edges, whose union is $V(H)$. 

As noted in \cite{vz}, not every hypergraph has a strongly minimal edge-cover (the set of all finite subsets of a countable set is a counterexample). In \cite{ak} it was shown that also not every hypergraph has a strongly maximal matching - a counterexample is the set of all sets of natural numbers whose cardinality and first element are the same. However, a basic problem is:

\begin{conjecture}\label{main}
In every hypergraph with  sizes of edges bounded from above by a natural number $k$ there exist a strongly maximal matching, a strongly minimal vertex-cover and a strongly minimal edge-cover. 
\end{conjecture}

The ``edge-cover'' version of the conjecture, together with a compactness argument, would imply the following conjecture of Aharoni and Korman, a result that is trivial in the finite case:

\begin{conjecture}[the fish bone conjecture]\label{fishbone}
Let $P$ be a poset with bounded width (namely having $|A|\le k$ for every antichain $A$, for some  finite number $k$)  Then there exist a chain $C$ and a decomposition $D$ of $V(P)$ into disjoint antichains, such that every antichain in $D$ meets $C$.

\end{conjecture}
In fact, we do not know a counterexample even just assuming having no infinite antichains.  \\

Conjectures \ref{main} and \ref{fishbone} are known for $k=2$ \cite{konig, tutte}. The difficulty of the proofs in this case indicates that the general conjectures are not easy.

Conjecture \ref{fishbone} is a special case of a more general conjecture, on infinite perfect graphs. Call a graph``perfect'' if all its induced finite subgraphs are perfect. Call a graph {\em strongly perfect} if in  every induced subgraph there exists a partition into independent sets, and a clique meeting them all. 

\begin{conjecture}\cite{al}
A perfect graph in which all independent sets are finite is strongly perfect.  
\end{conjecture}
In \cite{al} the conjecture was proved for chordal graphs and for their complements.

\end{document}